\documentclass[11pt]{amsart}
\usepackage{amsmath,amssymb,enumerate}
\usepackage[dvips]{color}
\date{October 27, 2009}   

 \numberwithin{equation}{section} \theoremstyle{plain}
\newcommand{\ds}{\displaystyle}
\newcommand{\tq}{\, \big| \, }
\renewcommand{\r}{\mathbb{R}}

\DeclareMathOperator{\vol}{vol}%
\DeclareMathOperator{\im}{im}%
\DeclareMathOperator{\comp}{comp}%
\DeclareMathOperator{\conf}{conf}%

 \theoremstyle{plain}    
\newtheorem{theorem}{\rm\bf Theorem}[section]
\newtheorem{proposition}[theorem]{\rm\bf Proposition}
\newtheorem{lemma}[theorem]{\rm\bf Lemma}
\newtheorem{corollary}[theorem]{\rm\bf Corollary}

\theoremstyle{definition}
\newtheorem{definition}{\rm\bf Definition}[section]
\makeatletter

\title{The H\"older-Poincar\'{e} Duality for $L_{q,p}$-cohomology}

\author{Vladimir Gol'dshtein} 
\address{Vladimir Gol'd'shtein, Department of Mathematics,  
Ben Gurion University of the Negev,  P.O.Box 653, Beer Sheva, Israel}  
\email{vladimir@bgu.ac.il}

\author{Marc Troyanov} 
\address{M. Troyanov, 
Section de MathŽmatiques,  
\'Ecole Polytechnique F{\'e}derale de Lausanne, 
1015 Lausanne - Switzerland}
\email{marc.troyanov@epfl.ch}

\begin{document}

\begin{abstract}
We prove the following version of Poincar\'{e} duality for reduced
$L_{q,p}$-cohomology: For any $1<q,p<\infty$ , the 
$L_{q,p}$-cohomology of a Riemannian manifold is in duality with the interior  $L_{p',q'}$\emph{-cohomology}
for $1/p+1/p'=1$, $1/q+1/q'=1$. This duality result is a generalization
of the corresponding result for $L_{p}$-cohomology \cite{GK4}. \medskip

\end{abstract}

\maketitle

\footnotetext{This work was partially supported by Israel Scientific Foundation (Grant No 1033/07)}

\section{Introduction and statement of the results}

The main goal of this paper is to describe the dual space of the reduced  $L_{q,p}$-{cohomology} of an oriented Riemannian manifold  $(M,g)$. 

Let us denote by \(\mathcal{D}^k = C^{\infty}_0(M,\Lambda^k)\) the vector space of smooth differential $k$-forms with compact support in $M$ and by $L^{p}(M,\Lambda^{k})$ the Banach space of $p$-integrable differential $k$-forms.
We also  consider the space  $ \Omega_{q,p}^{k}(M)$ of  $q$-integrable differential $k$-forms whose weak exterior differentials are $p$-integrable
$$
 \Omega_{q,p}^{k}(M) = \left\{ \,\omega\in L^{q}(M,\Lambda^{k})\,\big|\, d\omega\in L^{p}(M,\Lambda^{k+1})\right\}.
$$
 
We  now define our basic objects of investigation. 
\begin{definition}
The  reduced \emph{$L_{q,p}$-cohomology} of the Riemannian manifold $(M,g)$ is 
 defined as
\begin{equation}\label{eq.defRqp}
  \overline{H}^k_{q,p}(M) =  Z^k_{p}(M)/\overline{B}^k_{q,p}(M),
\end{equation}
where $Z^k_{p}(M)$ is the set of weakly closed forms in $L^{p}(M,\Lambda^{k})$
and  $\overline{B}^k_{q,p}$ is the closure in {$L^{p}(M,\Lambda^{k})$} of the set of weakly exact forms having 
a $q$-integrable primitive:
$$
 \overline{B}^k_{q,p} = \overline{d(\Omega_{q,p}^{k-1}(M))}^{L^{p}(M,\Lambda^{k})}.
 $$
\end{definition}
(we shall  use the notation $\overline{A}^{E}$ for the closure of the subset set $A$ in the Banach space $E$.)

\smallskip

In the special case $q = p$, the space defined in (\ref{eq.defRqp}) is simply
denoted by  $\overline{H}^k_{p}(M)=\overline{H}^k_{q,p}(M)$ and is called
the   reduced \emph{$L_{p}$-cohomology} of the manifold.

 \medskip

The reduced $L_{q,p}$-cohomology is naturally a Banach space. Two closed forms $\omega,\omega'$
in $Z_{p}^{k}(M)$ represent the same $L_{q,p}$-reduced cohomology
class if one can find a sequence $\theta_{j}\in L^{q}(M, \Lambda^{k-1}$ such that
$\Vert(\omega-\omega')-d\theta_{j}\Vert_{L^{p}(M,\Lambda^{k})\ }\rightarrow0$.

\bigskip 

The subject-matter of  $L_{q,p}$-cohomology has been the object of
a number of investigations in recent years. The paper  \cite{GT2006} contains some foundational material and shows how Sobolev inequalities for differential forms can be interpreted in the framework of $L_{q,p}$-cohomology , the paper  \cite{GT2009a} gives some applications to
quasi-conformal geometry  and \cite{GT2009b} relates the $L_{q,p}$-cohomology to more general classes of mappings. The paper  \cite{GT2009b} 
contains some computations for negatively curved Riemannian manifolds and the papers \cite{XLi2009,Troyanov2009} study the relation between
the $L_{q,p}$-cohomology of a manifold and the $L_p$-Hodge decomposition on that manifold.
In \cite{Kopylov2007,Kopylov2008}, some computations for warped product manifolds and the general Heisenberg group are developed.

\bigskip 

In order to describe the dual space to $\overline{H}^k_{q,p}(M)$, we introduce another
type of cohomology which we shall call the \emph{interior reduced $L_{r,s}$-cohomology}. 
This cohomology  captures the idea  of cohomology relative to the (ideal) boundary of the manifold. 
\begin{definition}
The  \emph{interior reduced $L_{r,s}$-cohomology} of the Riemannian manifold $(M,g)$ is the
Banach space defined as
\begin{equation}
  \overline{H}^k_{r,s;0}(M) =  Z^k_{r,s;0}(M)/\overline{d\mathcal{D}^{k-1}},  
\end{equation}
where $\overline{d\mathcal{D}^{k-1}}$ is the closure of $d\mathcal{D}^{k-1}$  in $L^{r}(M,\Lambda^{k})$ and 
$Z^k_{r,s;0}(M) \subset \Omega_{r,s}^{k}(M)$ is defined as:
$$
Z^k_{r,s;0}(M) =  \ker(d)\cap  \overline{\mathcal{D}^{k}(M)\ }^ 
{\Omega_{r,s}^{k}}.
$$
\end{definition}
A form $\alpha$ belongs to $Z^k_{r,s;0}(M)$ if and only if 
$\alpha$ is a weakly closed form in $L^{r}(M,\Lambda^{k})$ such that there exists a sequence 
$\theta_j \in \mathcal{D}^{k}(M)$ such that $\theta_j \to \alpha$ in $\Omega_{r,s}^{k}(M)$,
i.e.
$$
  \| \theta_j - \alpha\|_{r} \to 0 \qquad \mathrm{and} \qquad
  \| d\theta_j \|_{s} \to 0.
$$

Our main result is the following 

\begin{theorem}\label{th.duality1}
Let $(M,g)$ be an oriented $n$-dimensional Riemannian manifold. If $1<p,q<\infty$, then 
$\overline{H}^{k}_{q,p}(M)$ is isomorphic to the dual of  
$\,Ê\overline{H}^{n-k}_{p',q';0}(M)$
where $\frac{1}{p'}+\frac{1}{p}=\frac{1}{q'}+\frac{1}{q}=1$. The duality is induced by the
natural pairing $\overline{H}^{k}_{q,p}(M)\times \overline{H}^{n-k}_{p',q';0}(M)\rightarrow\r$
given by integration:
\begin{equation}\label{eq.chmpairing}
 \left\langle [\omega],[\varphi\right]\rangle =\int_{M}\omega\wedge\varphi.
\end{equation}
\end{theorem}

\medskip

It is useful to  introduce a third,  auxiliary,  cohomology which we called the \emph{semi-compact} or the
$(\comp,p)$-cohomology of $(M,g)$. This object  turns out to be more manageable than the 
reduced $L_{q,p}$-cohomology and its interior version.
\begin{definition} \label{def.compp}
The  \emph{reduced $(\comp, p)$-cohomology} of $(M,g)$ is the
Banach space defined as
$$
  \overline{H}^k_{\comp, p}(M) =  Z^k_{p}(M)/\overline{d\mathcal{D}^{k-1}},
$$
where $\overline{d\mathcal{D}^{k-1}}$ is the $L^p$ closure of $d\mathcal{D}^{k-1}$.
\end{definition}
Observe that the following inclusions hold
$$
  Z^k_{p,q;0}(M) \subset Z^k_{p}(M)
    \quad \text{and} \quad
  \overline{d\mathcal{D}^{k-1}} \subset \overline{B}^k_{q,p}.
$$
This implies that the reduced $L_{q,p}$-cohomology embeds in the  reduced $(\comp,p)$-cohomology and the latter
is a quotient of the reduced $L_{q,p}$-cohomology:
$$
   \overline{H}^k_{p,q;0}(M)  \hookrightarrow  \overline{H}^k_{\comp,p}(M) 
       \quad \text{and} \quad
     \overline{H}^k_{\comp,p}(M)  \twoheadrightarrow  \overline{H}^k_{q,p}(M) 
$$
It follows in particular that  $\dim ( \overline{H}^k_{q,p}(M) ) \leq \dim ( \overline{H}^k_{\comp,p}(M) )$
for any $q$.

\medskip

We  have  the following duality result for the $(\comp,p)$-cohomology.
\begin{theorem}\label{th.duality2}
Let $(M,g)$ be an oriented $n$-dimensional Riemannian manifold. If $1<p<\infty$, then 
$\overline{H}^{k}_{\comp,p}(M)$ is isomorphic to the dual of $\overline{H}^{n-k}_{\comp,p'}(M)$
where $\frac{1}{p'}+\frac{1}{p} = 1$. The duality is induced by the  integration pairing
(\ref{eq.chmpairing}).
\end{theorem}

\medskip

We now give a sufficient condition for the  reduced $(\comp,p)$-cohomology to coincide with the reduced $L_{q,p}$-cohomology. Recall that a Riemannian manifold  \((M,g)\) is said to be  \(s\)-\emph{parabolic}, \(1\leq s \leq \infty\) 
if one can approximate the unity by functions with small $s$-energy, i.e. 
 if there exists a sequence of smooth functions with compact support \(\{\eta_j\} \subset C^{\infty}_0(M)\) such that 
\(\eta_j \to 1\) uniformly on every compact subset of \(M\) and  \(\lim_{j\to \infty} \int_M |d\eta_j|^s = 0\).

\medskip

\begin{theorem}\label{th.parcomp}
Suppose that  \(M\) is $s$-parabolic for some $1\leq s \leq \infty$, and assume that
 \(\frac{1}{s} =  \frac{1}{p} - \frac{1}{q}\). Then we have 
 $$
  \overline{H}^k_{q,p}(M)   = \overline{H}^{k}_{\comp,p}(M)
       \quad \text{and} \quad
  \overline{H}^k_{q,p;0}(M) = \overline{H}^{k}_{\comp,q}(M).
$$      
\end{theorem}
This result  gives us in particular conditions under which the reduced $L_{q,p}$-cohomology and the 
interior $L_{q,p}$-cohomology coincide.

\medskip

The case $s= \infty$ is  important, because a manifold is  $\infty$-parabolic if and only if it is complete.
It follows that for complete manifolds, we have
$\overline{H}^k_{p}(M) = \overline{H}^k_{p;0}(M)  = \overline{H}^{k}_{\comp,p}(M)$ and   the pairing (\ref{eq.chmpairing}) induces a duality 
$$
  \overline{H}^k_{p}(M)' = \overline{H}^{n-k}_{p'}(M),
$$
where  $\frac{1}{p}+ \frac{1}{p'}  =  1$. The paper \cite{GT2009d} contains a short proof of this special case.

\medskip

Theorem \ref{th.parcomp} implies that if $M$ is $s$-parabolic for every $s$, then $\overline{H}^n_{q,p}(M)$ is independent of the choice of $q$ and  $\overline{H}^n_{q,p;0}(M)$ is independent of the choice of $p$ provided 
$p\leq q$. 

\medskip

In the finite volume case, we have a stronger result. If $(M,g)$ is complete with finite volume, then it is $s$-parabolic for
any $1\leq s \leq \infty$. For such manifolds, we then have
$$
   \overline{H}^k_{q,p}(M)  = \overline{H}^{k}_{\comp,p}(M) =   \overline{H}^k_{p}(M) 
$$
for any $q \geq p$. Combining this result with \cite[Proposition 7.1]{GT2006},
we obtain that $ \overline{H}^k_{q,p}(M) =0$ if and only if $\overline{H}^k_{p}(M)  =0$
(without the restriction $q \geq p$).

\bigskip

The rest of the paper is organized as follows: In the next section we recall some basic notions and facts and in  section 3,
we prove the duality Theorem \ref{th.duality2} for the $(\comp,p)$-cohomology. 
Section 4 is devoted to the proof of Theorem \ref{th.duality1} and it occupies the largest part of this paper. We reformulate
the problem in the convenient language of Banach complexes. This section also contains a description of the dual space
of  $\Omega^k_{q,p}(M)$, see Corollary \ref{cor.duality0}. In section 5 we recall the notion of $s$-parabolic manifolds
and we give a proof of Theorem \ref{th.parcomp}. In section 6 we deduce from the main theorem a result on the Poincar\'e duality
of the conformal cohomology and in  the last section, we derive some consequences and applications 
of our duality theorems.

\section{Some background}

For any  Riemannian $n$-manifold $(M,g)$, the Riemannian
metric induces a norm (in fact a scalar product) on $\Lambda^{k}T_{x}M^{*}$
at any point $x\in M$. Using this norm, one can define the space
$L^{p}(M,\Lambda^{k})$ of measurable differential $k$-forms $\omega$
such that
 \[
\Vert\omega\Vert_{p} = \left(\int_{M}|\omega|_{x}^{p}d\vol_{g}(x)\right)^{1/p}<\infty,\]
 if $1\leq p<\infty$ and  $\Vert\omega\Vert_{\infty} = \text{ess-sup} \,|\omega|<\infty$
 if $p = \infty$.
 If $p' = p/(p-1)$, one can define a pairing 
\[
\left\langle \,\,,\,\,\right\rangle :L^{p}(M,\Lambda^{k})\times L^{p'}(M,\Lambda^{n-k})\rightarrow \r
\]
 by integration, that is 
\begin{equation} \label{Ipairing}   
  \left\langle \omega,\varphi\right\rangle =\int_{M}\omega\wedge\varphi
\end{equation}
 for $\omega\in L^{p}(M,\Lambda^{k})$ and $\varphi\in L^{p'}(M,\Lambda^{n-k})$.
This is well defined because we have at (almost) every point $x\in M$
\[
\left|\omega\wedge\varphi\right|_{x}\leq\left|\omega\right|_{x}\left|\varphi\right|_{x},
\]
 and by H\"{o}lder's inequality 
 \[
\left|\left\langle \omega,\varphi\right\rangle \right| = 
\left|\int_{M}\omega\wedge\varphi\right| \leq
\int_{M}\left|\omega\right|\left|\varphi\right|d\mu
\leq 
\left\Vert \omega\right\Vert _{L^{p}(M,\Lambda^{k})}\left\Vert \varphi\right\Vert _{L^{p'}(M,\Lambda^{k})}.
\]

The bilinear function $\left\langle \,\,,\,\,\right\rangle $ allows
us to define for any $\varphi\in L^{p'}(M,\Lambda^{k})$ a bounded linear
functional  $F_{\varphi}(\omega) = \int_{M}\omega\wedge\varphi$, and the 
familiar H\"older duality between $p$ and $p'$ extends to differential forms,
see e.g. \cite{GK4}:
\begin{theorem} \label{thm:Lduality}
For any \(1\leq p<\infty\) the correspondence
\(\varphi\to F_{\varphi}\) is an isometric  isomorphism from the Banach
spaces \(L^{p'}(M,\Lambda^{n-k})\) to the dual space of \(L^{p}(M,\Lambda^{k})\),
and the duality is explicitly given by the pairing (\ref{Ipairing}).
\end{theorem}

\bigskip
 
\begin{corollary}\label{cor.density1}
\(\mathcal{D}^k := C^{\infty}_0(M,\Lambda^k)\) is dense in  \(L^{p}(M,\Lambda^k)\).
\end{corollary}

\textbf{Proof.}  
Here is a short proof: Suppose there exists a $k$-form $\omega\in L^{p}(M,\Lambda^{k})$
such that $\omega\not\in\overline{\mathcal{D}}^{k}$. By the Hahn-Banach Theorem,
there exists a continuous linear form $\lambda:L^{p}(M,\Lambda^{k})\rightarrow\r{}$
such that $\lambda=0$ on $\overline{\mathcal{D}}^{k}$ and $\lambda(\omega)\neq0$.
By the previous Theorem, there exists $\psi\in L^{p'}(M,\Lambda^{n-k})$
such that $\lambda(\theta)=\int_{M}\theta\wedge\psi$ for any $\theta\in L^{p}(M,\Lambda^{k})$.
In particular $\lambda(\theta)=\int_{M}\theta\wedge\psi=0$ for any
$\theta\in\mathcal{D}^{k}$. This implies that $\psi=0$ and therefore
$\lambda(\omega)=0$. Therefore  no such $\omega$ exists and
we conclude that $\overline{\mathcal{D}}^{k}=L^{p}(M,\Lambda^{k})$. 

\qed

\bigskip

We now define the notion of weak exterior differential:
\begin{definition}\label{def.weakd}
Assume $M$ to be oriented. One says that a differential form  $\theta \in L^p(M,\Lambda^{k+1})$ is the
\emph{weak exterior differential} of the form $\phi \in L^p(M,\Lambda^{k})$ and one writes
{$d\phi = \theta$} if one has
\[
 \int_M \theta \wedge \omega = (-1)^{k+1}\int_M \phi \wedge d\omega
\]
for any  $\omega \in\mathcal{D}^{n-k}(M)= C^{\infty}_0(M,\Lambda^{n-k})$.
\end{definition}

We then introduce the space 
\[
\Omega_{q,p}^{k}(M) = \left\{ \,\omega\in L^{q}(M,\Lambda^{k})\,\big|\, d\omega\in L^{p}(M,\Lambda^{k+1})\right\} ,
\]
 this is a Banach space for the graph norm 
 \begin{equation}\label{grafnorm}
\left\Vert \omega\right\Vert _{q,p} = \left\Vert \omega\right\Vert _{L^{q}}+\left\Vert d\omega\right\Vert _{L^{p}}.
\end{equation}
The space $\Omega_{q,p}^{k}(M)$ is a reflexive Banach
space for any $1<q,p<\infty$, this can be proved using  standard arguments of functional analysis (see e.g. \cite{Brezis}),
but is  also follows from Corollary \ref{cor.duality0} below.

 \medskip

\begin{lemma}\label{lem.density1}
For any  \(1\leq q,p<\infty\),  \(C^{\infty}(M,\Lambda^{k})\cap\Omega_{q,p}^{k}(M)\) is dense in \(\Omega_{q,p}^{k}(M)\).
\end{lemma}

\textbf{Proof.}
This follows from the regularization theorem, see \cite{GT2006}. 

\qed

\bigskip

We now (re)define our basic ingredients, for $p,q,r \in [1,\infty]$.
\begin{definition}  The closure of $\mathcal{D}^{k}=C_{0}^{\infty}(M,\Lambda^{k})$ in $\Omega_{q,p}^{k}(M)$
is denoted by $\Omega_{q,p;0}^{k}(M)$. We also define the following subspaces:
\begin{enumerate}[(a)]
 \medskip \item $Z^k_{p,r}(M) = \ker [d : \Omega_{p,r}^{k}(M) \to L^{r}(M,\Lambda^{k+1}) ]$.
 \medskip  \item $B^k_{q,p}(M) = \im [d : \Omega_{q,p}^{k-1}(M)\to L^{p}(M,\Lambda^{k}) ]$.
 \medskip \item $Z^k_{p,r;0}(M) = \ker [d : \Omega_{p,r;0}^{k}(M)\to L^{r}(M,\Lambda^{k+1}) ]$.
 \medskip  \item $B^k_{q,p;0}(M) = \im [d : \Omega_{q,p;0}^{k-1}(M)\to L^{p}(M,\Lambda^{k}) ]$. 
  \end{enumerate}
\end{definition}

\medskip

\begin{lemma}\label{lem.density2}
The previously defined spaces satisfy the following properties
\begin{enumerate}[(i)]
  \item $Z^k_{p,r}(M)$ does not depend on $r$ and it is is a closed subspace of $L^{p}(M,\Lambda^{k})$.
  \item $d(\mathcal{D}^{k-1})$ is dense in $B_{r,p;0}^{k}(M)$ for the $L^p$-topology.
\end{enumerate}
\end{lemma}

\textbf{Proof.}
(i) $Z^k_{p,r}(M)$ is a closed subspace of $\Omega_{p,r}^k(M)$ since it is the kernel of the bounded operator $d$.
It is also a closed subspace of $L^{p}(M,\Lambda^{k})$ since 
 for any $\alpha \in Z^k_{p,r}(M) $, we have  $||\alpha||_{\Omega_{p,r}^k(M)}= ||\alpha||_{L^{p}(M,\Lambda^{k})}$.
 Now  $Z^k_{p,r}(M)$ does not depend on $r$ because it coincides with the space of
weakly closed $k$-forms in $L^{p}(M,\Lambda^{k})$.

\medskip

Statement (ii) is  almost obvious, let $\omega\in B_{r,p;0}^{k}(M)$. Then there
exists $\theta\in\Omega_{r,p;0}^{k-1}(M)$ such that $d\theta=\omega$.
By definition of $\Omega_{r,p;0}^{k-1}(M)$, there exists a sequence
$\theta_{j}\in \mathcal{D}^{k-1}$ such that 
$\lim_{j\rightarrow\infty}\Vert\theta-\theta_{j}\Vert_{r,p}=0$.
This means that $\theta=\lim_{j\rightarrow\infty}\theta_{j}$ in $L^{r}$
and $\omega=d\theta=\lim_{j\rightarrow\infty}d\theta_{j}$ in $L^{p}$. 

\qed

\bigskip

The Banach space $Z^k_{p,r}(M)$ will then simply be denoted by $Z^k_{p}(M)$. We will also identify the closed
subspaces $\overline{B}_{r,p;0}^{k}(M)$ and $\overline{d(\mathcal{D}^{k-1})}$. 
Our reduced cohomologies are then naturally defined as the following quotients of Banach spaces:
$$
 \overline{H}^k_{q,p}(M) = Z^k_{p}(M)/\overline{B}^k_{q,p}(M), \qquad
 \overline{H}^k_{p,q;0}(M) = Z^k_{p,q;0}(M)/\overline{d(\mathcal{D}^{k-1})}(M)
$$
and
$$
 \overline{H}^k_{\comp,p}(M) = Z^k_{p}(M)/\overline{d(\mathcal{D}^{k-1})}.
$$
We also defined the \emph{unreduced} $L_{q,p}$-cohomology as
$$
  H^k_{q,p}(M) = Z^k_{p}(M)/B^k_{q,p}(M),  
$$
this is generally not a Banach space as $B^k_{q,p}(M) \subset Z^k_{p}(M)$ need not be a closed subspace.
We may also define an unreduced interior cohomolgy, but this space will depend on 3 indices instead of 2:
$$
 H^k_{r,p,q;0}(M) = Z^k_{p,q;0}(M)/B_{r,p;0}^{k}(M).
$$

\section{The duality theorem for the $(\comp,p)$-cohomology}

In this section, we prove Theorem \ref{th.duality2}.
The proof is based on the following lemma from functional analysis.

\medskip

\begin{lemma}\label{lem.dualquotient}
Let \(I:X_{0}\times X_{1}\to \mathbb{R}\) be a duality (non degenerate pairing) between two reflexive Banach spaces. Let 
\(B_{0},A_{0},B_{1},A_{1}\) be linear subspaces  
such that
\[
 B_{0}\subset A_{0} = B_{1}^{\bot} \subset X_{0}
       \qquad  \text{and} \qquad  
B_{1}\subset A_{1} = B_{0}^{\bot}\subset X_{1}.
\]
Then the pairing \(\overline{I}:\overline{H}_{0}\times\overline{H}_{1}\to \mathbb{R}\)
of \(\overline{H}_{0} = A_{0}/\overline{B}_{0}\) and \(\overline{H}_{1} = A_{1}/\overline{B}_{1}\)
is well defined and induces duality between \(\overline{H}_{0}\) and
\(\overline{H}_{1}\).
\end{lemma}
Here $B_{1}^{\bot}$ and $B_{0}^{\bot}$ are the annihilators of $B_1$ and $B_0$ for the duality $I$.

\bigskip

\textbf{Proof.}  We first prove the Lemma under the assumption that
$A_0  \neq \overline{B}_0$ and  $A_1  \neq \overline{B}_1$. 
Observe  that $A_i$ is a closed subspace of the Banach spaces $X_i$ since it 
is the annihilators of $B_i$ ($i = 1,2$). 
The bounded bilinear map $I:A_{0}\times A_{1}\rightarrow\mathbb{R}$
is defined by restriction and it gives rise to a well defined bounded
bilinear map 
$$\overline{I}:A_{0}/\overline{B}_{0}\times A_{1}/\overline{B}_{1}\rightarrow\mathbb{R}$$
because we have the inclusions $B_{0}\subset A_{0}\subset B_{1}^{\bot}$
and $B_{1}\subset A_{1}\subset B_{0}^{\bot}$.
(Indeed, let $\alpha_0 \in A_0, \alpha_1 \in A_1, b_0 \in \overline{B_0}, b_1 \in \overline{B_1}$. Then 
$$
I(\alpha_0+b_0, \alpha_1+b_1)=I(\alpha_0,\alpha_1)+I(\alpha+b_0,b_1)+I(b_0,\alpha_1),
$$
by the definition of annihilators the second and third terms vanish and the duality
$\overline{I}$ is thus well defined.)

\smallskip

We show that $\overline{I}$ is non degenerate: let $a_{0}\in A_{0}$
be such that $[a_0]\neq0\in A_{0}/\overline{B}_{0}$; i.e. $a_0\not\in\overline{B}_{0}$.
By Hahn-Banach theorem and the fact that $X_{1}$ is dual to $X_{0}$, there
exists an element $y\in X_{1}$ such that $I(a_0,y)\neq0$ and $I(b,y)=0$
for all $b\in\overline{B}_{0}$. Thus $y\in B_{0}^{\bot}=A_{1}$ and
we have found an element $[y]\in A_{1}/\overline{B}_{1}$ such that
$\overline{I}([a],[y])\neq0$.

The same argument shows that for any $[\alpha]\neq0\in A_{1}/\overline{B}_{1}$,
we can find an element $[x]\in A_{0}/\overline{B}_{0}$ such that
$\overline{I}([x],[\alpha])\neq0$. The proof is thus complete in the case
 $A_i  \neq \overline{B}_i$. 
 
 \smallskip
 
If $A_0 = \overline{B}_0$, then we have from the hypothesis of the Lemma:
$$
 \overline{B}_{1} \subset A_{1} = B_{0}^{\bot} = A_{0}^{\bot} 
 = (B_{1}^{\bot})^{\bot} =  \overline{B}_{1},
$$
and it follows that $A_1 = \overline{B}_1$.  The argument shows in fact  that
$$
 A_0 = \overline{B}_0  \quad  \Leftrightarrow \quad
 A_1 = \overline{B}_1.
$$
In that case both  \(\overline{H}_{0} = A_{0}/\overline{B}_{0}\) and \(\overline{H}_{1} = A_{1}/\overline{B}_{1}\)
are null spaces and the Lemma is trivial.

\qed

\bigskip
 
\textbf{Proof of Theorem \ref{th.duality2}}
Let $\alpha \in L^{p}(M,\Lambda^{k})$, then by definition of the weak exterior differential, we have 
$d\alpha = 0$ if and only if
$$
 \int_M \alpha \wedge d\omega = 0
$$
for any $\omega \in \mathcal{D}^{n-k-1}$. This precisely means that $Z^k_{p}(M) \subset  L^{p}(M,\Lambda^{k})$
is the annihilator of ${d\mathcal{D}^{n-k-1}}$ for the integration pairing, we thus have the following relations with respect to  the pairing  (\ref{Ipairing}):
\[
    dD^{k-1}\subset Z^k_{p}(M) =  ({d\mathcal{D}}^{n-k-1})^{\perp} \subset  L^{p}(M,\Lambda^{k}).
\]
Likewise, we have 
$\ dD^{n-k-1}\subset  Z^{n-k}_{p'}(M) =  ({d\mathcal{D}}^{k-1})^{\perp} \subset L^{p'}(M,\Lambda^{n-k})$, \ 
and the previous Lemma implies that  (\ref{Ipairing}) induces a duality between  \\
$\overline{H}^k_{\comp,p}(M) =  Z^k_{p}(M)/\overline{d\mathcal{D}^{k-1}}$ and
 $\overline{H}^{n-k}_{\comp,p}(M) =  Z^k_{p}(M)/\overline{d\mathcal{D}^{n-k-1}}$.
 
 \qed

\section{The main duality Theorem}

The goal of this section  is to prove Theorem \ref{th.duality1}. 
It will be convenient to use the language of Banach complexes, and for that we need to reformulate the definition of our cohomology in a new language.

\subsection{The $L_{\pi}$-cohomology of $M$}

To define a Banach complex, we fix an $(n+1)$-tuple of real numbers
\begin{equation}\label{eq.pi}
\pi=\{ p_{0},p_{1},\cdots,p_{n}\}\subset[1,\infty]
\end{equation}
 and define 
 \[
\Omega_{\pi}^{k}(M) = \Omega_{p_{k},p_{k+1}}^{k}(M).
\]
Observe that $\Omega_{\pi}^{n}(M)=L^{p_{n}}(M,\Lambda^{n})$ and $\Omega_{p,p}^{1}(M)$
coincides with the Sobolev space $W^{1,p}(M)$.
Since the weak exterior differential is a bounded operator \
$d:\Omega_{\pi}^{k-1}\rightarrow\Omega_{\pi}^{k}$, we have constructed
a Banach complex $\{\Omega_{\pi}^{*}(M),d\}$: 
\[
0\rightarrow\Omega_{\pi}^{0}\overset{d}{\rightarrow}\cdots\overset{d}{\rightarrow}\,\,
\Omega_{\pi}^{k-1}\overset{d}{\rightarrow}\,\,\Omega_{\pi}^{k}\overset{d}{\rightarrow}\cdots\overset{d}{\rightarrow}\
\Omega_{\pi}^{n}\rightarrow0\,.
\]

\begin{definition}
\emph{The (reduced)} \(L_{\pi}\)-\emph{cohomology of} \(M\) \emph{is
the (reduced) cohomology of the Banach complex} \(\{\Omega_{\pi}^{k}(M),d_{k}\}\).
\end{definition} 

The $L_{\pi}$-cohomology space $H_{\pi}^{k}(M)$ depends only on
$p_{k\text{ }}$and $p_{k-1}$ and we have in fact 
$$
 H_{\pi}^{k}(M) = H_{p_{k-1},p_{k}}^{k}(M) =Z_{p_{k}}^{k}(M)/B_{p_{k-1},p_{k}}^{k}(M)
$$ 
and  the reduced $L_{\pi}$-cohomology is
$$
\overline{H}_{\pi}^{k}(M)= \overline{H}_{p_{k-1},p_{k}}^{k}(M)
 =Z_{p_{k}}^{k}(M)/\overline{B}_{p_{k-1},p_{k}}^{k}(M).
$$

\medskip

We also introduce a notion of interior $L_{\pi}$-cohomology.
Let us denote by $\Omega_{\pi;0}^{k}(M)$ the closure of $\mathcal{D}^k(M)$
in $\Omega_{\pi}^{k}(M)$. This is  an another Banach complex 
\[
0\rightarrow\Omega_{\pi;0}^{0}(M)\overset{d}{\rightarrow}\cdots\overset{d}{\rightarrow}\,\,
\Omega_{\pi;0}^{k-1}(M)\overset{d}{\rightarrow}\,\,
\Omega_{\pi;0}^{k}(M)\overset{d}{\rightarrow}\cdots\overset{d}{\rightarrow}\
\Omega_{\pi;0}^{n}\rightarrow0\,.
\]

\begin{definition}
The \emph{interior \(L_{\pi}\)-cohomology } of \((M,g)\)  is the cohomology of this
Banach complex, i.e.  
 \[
H_{\pi;0}^{k}(M)=H_{p_{k-1},p_{k},p_{k+1};0}^{k}(M) = Z_{p_{k},p{k+1}:0}^{k}(M)/B_{p_{k-1},p_{k}:0}^{k}(M)
\]
 and by Lemma \ref{lem.density2} the interior reduced $L_{\pi}$-cohomology  is
\[
\overline{H}_{\pi;0}^{k}(M)=\overline{H}_{p_{k},p_{k+1};0}^{k}(M)
 = Z_{p_{k},p{k+1}:0}^{k}(M)/\overline{dD^{k-1}}(M).
\]
\end{definition}

\medskip

\begin{definition}
The \emph{dual} of the  $(n+1)$-tuple $\pi=\{ p_{0},p_{1},...,p_{n}\}\subset[1,\infty]$ 
 Is the $(n+1)$-tuple $\pi'=\{ r_{0},r_{1}...,r_{n}\}$ such that $\frac{1}{r_{k}}+\frac{1}{p_{n-k}}=1$.
\end{definition}
In the sequel, we will  use  the notation $a<\pi<b$ or $a\leq\pi\leq b$ if these inequalities hold for all $p_{0},p_{1},\cdots,p_{n}$.

\subsection{The complex of pairs of forms on $M$} 

Because the space $\Omega_{\pi}^{k}(M)$ is equipped with the graph
norm, it is useful to investigate the structure of the space where
this graph leaves, that is the following Cartesian product:

\begin{definition}
Given a $n$-dimensional Riemannian manifold $(M,g)$ and a sequence 
$\pi$ as in (\ref{eq.pi}), we introduce the
vector space 
\[
\mathcal{P}_{\pi}^{\star}(M) = L_{\pi}(M,\Lambda^{\star})\oplus L_{\pi}(M,\Lambda^{\star}).
\]
\end{definition}
We introduce a graduation of this vector space  be defining 
\[
\mathcal{P}_{\pi}^{k}(M) = L_{\pi}(M,\Lambda^{k})\oplus L_{\pi}(M,\Lambda^{k+1}).
\]

\medskip

\textbf{Remark} \ We have $L_{\pi}(M,\Lambda^{n+1})=L_{\pi}(M,\Lambda^{-1})=\{0\}$,
hence $\mathcal{P}_{\pi}^{n}(M)\cong L^{p_{n}}(M,\Lambda^{n})$ and
$\mathcal{P}_{\pi}^{-1}(M)\simeq L^{p_{0}}(M,\Lambda^{0})$.

\medskip

An element in $\mathcal{P}_{\pi}^{k}(M)$ will be denoted as a column
vector $\binom{\alpha}{\beta}$ with $\alpha\in L_{\pi}(M,\Lambda^{k})$
and $\beta\in L_{\pi}(M,\Lambda^{k+1})$. This space is a Banach space
for the norm 
\begin{equation}
\left\Vert \binom{\alpha}{\beta}\right\Vert _{\mathcal{P}_{\pi}^{k}(M)} = 
\left\Vert \alpha\right\Vert _{p_{k}}+\left\Vert \beta\right\Vert _{p_{k+1}},\label{norminP}
\end{equation}
 and it can be turned as a Banach complex for the ``differential''
$d_{\mathcal{P}}:\mathcal{P}_{\pi}^{k}\rightarrow\mathcal{P}_{\pi}^{k+1}(M)$
defined by 
\[
d_{\mathcal{P}}\binom{\alpha}{\beta}=\binom{\beta}{0}.
\]

\medskip

\begin{lemma}\label{lem:Ptrivial}
 The complex \((\mathcal{P}^{k}_{\pi}(M) , d _{\mathcal{P}})\) has trivial cohomology.
\end{lemma}

\textbf{Proof} Suppose that $d_{\mathcal{P}}\binom{\alpha}{\beta}=\binom{0}{0}$,
this means that $\beta=0$. But then it is clear that 
\[
 \binom{\alpha}{\beta}=\binom{\alpha}{0} =d_{\mathcal{P}}\binom{0}{\alpha}.
\]
\qed

\bigskip

\begin{proposition}\label{pro:LLduality} 
Let \(\pi\) be a sequence as in the Definition (\ref{eq.pi}) and \(\pi'\) be the dual sequence. 
If $M$ is oriented and $1\leq \pi < \infty$, then \(\mathcal{P}^{k}_{\pi}(M)\)
and \(\mathcal{P}^{n-k-1}_{\pi'}(M)\) are in duality for the following pairing 
\begin{equation} \label{eq:pairing}
\left\langle \binom{\alpha}{\beta} , \binom{\varphi}{\psi}  \right\rangle =
\int_M ((-1)^{k}\alpha\wedge \psi + \beta\wedge \varphi).
\end{equation}
\end{proposition}

\textbf{Proof.}
By Theorem \ref{thm:Lduality}, we have the following isomorphims
(in fact isometries): 
\[
\left(L^{p_{k}}(M,\Lambda^{k})\right)^{\prime}\cong L^{p'_{k}}(M,\Lambda^{n-k}),  \quad
\left(L^{p_{k+1}}(M,\Lambda^{k+1})\right)^{\prime}\cong L^{p'_{k+1}}(M,\Lambda^{n-k-1})
\]
 induced by isometrical isomorphisms. Hence 
\begin{align*}
\left(\mathcal{P}_{\pi}^{k}(M)\right)^{\prime}
 & =\left(L^{p_{k}}(M,\Lambda^{k})\oplus L^{p_{k+1}}(M,\Lambda^{k+1})\right)^{\prime}\\
 & =\left(L^{p_{k}}(M,\Lambda^{k})\right)^{\prime}\oplus\left(L^{p_{k+1}}(M,\Lambda^{k+1})\right)^{\prime}\\
 & =L^{p'_{k+1}}(M,\Lambda^{n-k-1})\oplus L^{p'_{k}}(M,\Lambda^{n-k})\\
 & =\mathcal{P}_{\pi'}^{n-k-1}(M).
 \end{align*}
\qed

\begin{lemma}\label{lemadj1}
The  operator  \(d _{\mathcal{P}}: \mathcal{P}_{\pi'}^{k-1} \to \mathcal{P}_{\pi'}^{k}\)
  and  \((-1)^{k}d _{\mathcal{P}}: \mathcal{P}_{\pi}^{n-k-1} \to \mathcal{P}_{\pi}^{n-k}\)
  are adjoint for the duality  (\ref{eq:pairing}).
\end{lemma}

\textbf{Proof.}
Let $ \binom{\alpha}{\beta} \in  \mathcal{P}_{\pi'}^{k-1}$ and 
$\binom{\varphi}{\psi} \in  \mathcal{P}_{\pi}^{n-k-1}$, then 
$$
 \left\langle d _{\mathcal{P}}\binom{\alpha}{\beta} , \binom{\varphi}{\psi}  \right\rangle 
 =  \left\langle \binom{\beta}{0} , \binom{\varphi}{\psi}  \right\rangle 
 = \int_M (-1)^{k}\beta\wedge \psi .
$$
On the other hand
$$
 \left\langle \binom{\alpha}{\beta} , d _{\mathcal{P}}\binom{\varphi}{\psi}  \right\rangle 
 =  \left\langle \binom{\alpha}{\beta}  , \binom{\psi}{0}   \right\rangle 
 = \int_M \beta\wedge \psi .
$$
\qed

\subsection{The annihilator of $\Omega_{\pi;0}^{\star}(M)$ in $\mathcal{P}_{\pi}^{\star}(M)$ }

We will investigate  $\Omega_{\pi}^{\star}(M)$ as a subspace of $\mathcal{P}_{\pi}^{\star}(M)$.
\begin{definition}
We denote by ${\Sigma}_{\pi}^{\star}(M) \subset \mathcal{P}_{\pi}^{\star}(M)$
the set of pairs of the form $\binom{\omega}{d\omega}$, where $\omega \in \Omega_{\pi}^{\star}(M)$.
We also denote by${\Sigma}_{\pi;0}^{\star}(M)$ the suspace of those elements $\binom{\omega}{d\omega}$
such that  $\omega \in \Omega_{\pi;0}^{\star}(M)$.
\end{definition}

It is clear that $\Sigma_{\pi}^{\star}(M)$ and $\Sigma_{\pi;0}^{\star}(M)$ are closed subspaces
of $\mathcal{P}_{\pi}^{\star}(M)$, they are subcomplexes and are isomorphic (as Banach complexes)
to $\Omega_{\pi;0}^{\star}(M)$ and $\Omega_{\pi}^{\star}(M)$.

\begin{lemma}
  The subspaces \(\Sigma^{k}_{\pi}(M) \subset \mathcal{P}^{k}_{\pi}(M)\) 
  and \(\Sigma^{n-k-1}_{\pi';0}(M) \subset \mathcal{P}^{n-k-1}_{\pi'}(M)\) 
  are orthogonal for the duality pairing (\ref{eq:pairing}).
\end{lemma}

\medskip

\textbf{Proof.}
 This is clear by the definition of weak exterior differential
and the density of $\mathcal{D}^{n-k-1}$ in $\Omega_{\pi';0}^{n-k-1}(M)$.

\qed

\medskip

In fact we have a stronger result:

\begin{proposition}
The subspace \(\Sigma^{k}_{\pi}(M)\) in \(\mathcal{P}^{k}_{\pi}(M)\) is the 
annihilator of \(\Sigma^{n-k-1}_{\pi;0}(M)\):
\begin{equation} \label{eq.ann1}
 \Sigma^{k}_{\pi}(M) = (\Sigma^{n-k-1}_{\pi';0}(M))^{\perp}.
\end{equation}
We also have
\begin{equation} \label{eq.ann2}
  \Sigma^{k}_{\pi; 0}(M) = (\Sigma^{n-k-1}_{\pi'}(M))^{\perp}.
\end{equation}
\end{proposition}

\medskip

\textbf{Proof} By density of  $\mathcal{D}^{n-k-1}$  in
$\Omega_{\pi;0}^{n-k-1}(M)$,  an element $\binom{\alpha}{\beta}\in\mathcal{P}_{\pi}^{k}(M)$
belongs to the annihilator of $\mathcal{D}^{n-k-1}$ if and
only if 
\[
(-1)^{k}\int_{M}\alpha\wedge d\varphi+\int_{M}\beta\wedge\varphi=0
\]
for any $\mathcal{D}^{n-k-1}$, but this means
by definition that $\beta=d\alpha$ in the weak sense, i.e. that $\binom{\alpha}{\beta}\in\Sigma_{\pi}^{k}(M)$.
Therefore 
\[
\Sigma_{\pi}^{k}(M)=\left(\mathcal{D}^{n-k-1}\right)^{\bot}=\left(\Sigma_{\pi';0}^{n-k-1}(M)\right)^{\bot},
\]
 this proves (\ref{eq.ann1}). To prove (\ref{eq.ann2}), we use the
fact that $\Sigma_{\pi';0}^{n-k-1}(M)\subset\mathcal{P}_{\pi'}^{n-k-1}(M)$
is a closed subspace  together with (\ref{eq.ann1})
and the following property of annihilator $(A^{\bot})^{\bot} =\overline{A}$
to deduce that 
\[
\Sigma_{\pi;0}^{k}(M)=\left(\left(\Sigma_{\pi;0}^{k}(M)\right)^{\bot}\right)^{\bot}=\left(\Sigma_{\pi'}^{n-k}(M)\right)^{\bot}.
\]
\qed

\subsection{The dual of $\Omega_{\pi}^{\star}(M)$ }
We introduce the following quotient of $\mathcal{P}_{\pi'}^{\star}(M)$:
\begin{equation}
\mathcal{A}_{\pi'}^{\star} = \mathcal{P}_{\pi'}^{\star}(M)/\Sigma_{\pi';0}^{\star}(M),\label{eq.defA}.
\end{equation}
This space inherits a graduation from that of   $\mathcal{P}_{\pi'}^{\star}(M)$.

\medskip

\begin{theorem}  \label{th.propertyofA}
The graded vector space \(\mathcal{A}_{\pi'}^{\star}\) has the following properties:
\begin{enumerate}[a.)] 
\item \(\mathcal{A}_{\pi'}^{k}\) is a Banach space for the norm
\[
\left\Vert \binom{\varphi}{\psi}\right\Vert_{ \mathcal{A}} =
\inf\left\{ \left(\left\Vert \varphi+\rho \right\Vert _{p'_{k}}+\left\Vert \psi +d\rho\right\Vert _{p'_{k+1}}\right)
\tq \rho \in\Omega_{\pi';0}^{k}(M)\right\}.
\]
\smallskip
\item \(\mathcal{A}_{\pi'}^{k}\)  is dual to \  \({\Sigma}_{\pi}^{n-k-1}(M)\) for the pairing given by (\ref{eq:pairing}).
\smallskip
\item The differential Ê \(d_{\mathcal{P}} : \mathcal{P}_{\pi'}^{k} \to \mathcal{P}_{\pi'}^{k+1}\) induces a differential  \(d_{\mathcal{A}} : \mathcal{A}_{\pi'}^{k} \to \mathcal{A}_{\pi'}^{k+1}\) and 
\(\left(\mathcal{A}_{\pi'}^{\star} , d_{\mathcal{A}}\right)\) is a Banach complex.
\smallskip
\item The  operator  \(d_{\mathcal{A}} : \mathcal{A}_{\pi'}^{k-1} \to \mathcal{A}_{\pi'}^{k}\)
  and  \(d : {\Sigma}_{\pi}^{n-k-1} \to {\Sigma}_{\pi}^{n-k}\)
  are adjoint (up to sign) for the duality  (\ref{eq:pairing}).
\end{enumerate}
\end{theorem}

\textbf{Proof.}
 The statement (a) is a standard fact on quotient of Banach spaces and (b) is a consequence of
 the orthogonality relation (\ref{eq.ann2}).
 The statement (c) follows directly from the fact that ${\Sigma}_{\pi';0}^{*}$ is a subcomplex of
 $\mathcal{P}_{\pi'}^{*}$, i.e.  $d_{\mathcal{P}}({\Sigma}_{\pi';0}^{k} )
 \subset ({\Sigma}_{\pi';0}^{k+1})$. Finally, the last statement follows directly from Lemma \ref{lemadj1}.
 
\qed

\medskip

\begin{corollary}\label{cor.duality0}
 The dual of the space \({\Omega}_{\pi}^{k}(M)\) is isomorphic to the completion of $\mathcal{D}^{n-k}$ with respect to
 the following norm:
\[
\left\Vert \sigma   \right\Vert =
\inf\left\{ \left(\left\Vert \rho \right\Vert _{p'_{n-k-1}}+\left\Vert \sigma +d\rho\right\Vert _{p'_{n-k}}\right)
\tq \rho \in\Omega_{\pi';0}^{n-k-1}(M)\right\}.
\]
\end{corollary}

\textbf{Proof.}
We know be the  Corollary \ref{cor.density1} that the space 
$\mathcal{S}^m = \mathcal{D}^{m} \oplus \mathcal{D}^{m+1}$ is a dense subspace of $\mathcal{P}^{m}$.
The image of $\mathcal{S}^m$ in $\mathcal{A}_{\pi'}^{m}$ is thus also dense. Observe that $\binom{\varphi}{\psi} \in \mathcal{S}^m$ and $\binom{0}{\psi - d\varphi} = \binom{\varphi}{\psi}  - \binom{\varphi}{d\varphi} \in \mathcal{S}^m$ represent the same element in $\mathcal{A}_{\pi'}^{m}$; this implies that the map 
$$
  j : \mathcal{D}^{m+1}  \to \mathcal{A}_{\pi'}^{m}, \quad 
    j(\sigma) =   \binom{0}{\sigma}
$$
has a dense image. This map is furthermore injective, indeed, if $\sigma \in \ker j$, then $\binom{0}{\sigma}
=  \binom{\rho}{\rho}$ for some $\rho \in \Sigma_{\pi';0}^{m-1}(M)$, but then $\rho = 0$ and therefore 
$\sigma = d\rho = 0$. It follows that $\mathcal{A}_{\pi'}^{m}$ is the completion of $j(\mathcal{D}^{m+1})$
for the natural norm in $\mathcal{A}_{\pi'}^{m}$. This norm is given by 
\[
\left\Vert j(\sigma)   \right\Vert 
= \left\Vert \binom{0}{\sigma}  \right\Vert_{\mathcal{A}} = 
\inf\left\{ \left(\left\Vert \rho \right\Vert _{p'_{m-1}}+\left\Vert \sigma +d\rho\right\Vert _{p'_{m}}\right)
\tq \rho \in\Omega_{\pi';0}^{m-1}(M)\right\}.
\]
We now deduce from statement (b) of the previous Theorem  that the above completion of $\mathcal{D}^{m+1}$
is isomorphic to the dual of ${\Omega}_{\pi}^{n-m-1}(M)$. Setting $m = n-k-1$ completes the proof of the Corollary.

\qed

\medskip

\subsection{The  duality theorem in $L_{\pi}$-cohomology}

In this subsection, we prove a duality statement between the reduced $L_{\pi}$-cohomology
of $M$ and the reduced interior $L_{\pi'}$-cohomology. We begin by investigating the cohomology of the 
 Banach complex  $\mathcal{A}_{\pi'}^{\star}$. Let us define
 the spaces
 $$
  Z^k(\mathcal{A}_{\pi'}^{\star}) = \ker[d_{\mathcal{A}} : \mathcal{A}_{\pi'}^{k} \to \mathcal{A}_{\pi'}^{k+1}]
\quad \text{and} \quad
  B^k(\mathcal{A}_{\pi'}^{\star}) = \im[d_{\mathcal{A}} (\mathcal{A}_{\pi'}^{k-1})].
 $$
We also denote by $\overline{B}^k(\mathcal{A}_{\pi'}^{\star})$ the closure of $B^k(\mathcal{A}_{\pi'}^{\star})$. 
Observe that $Z^k(\mathcal{A}_{\pi'}^{\star}) \subset \mathcal{A}_{\pi'}^{\star}$ is a closed subspace
and that  $B^k(\mathcal{A}_{\pi'}^{\star})  \subset \overline{B}^k(\mathcal{A}_{\pi'}^{\star})  \subset 
Z^k(\mathcal{A}_{\pi'}^{\star})$. The cohomology and reduced cohomology of  $\mathcal{A}_{\pi'}^{\star}$
are defined as
$$
 H^k(\mathcal{A}_{\pi'}^{\star}) = Z^k(\mathcal{A}_{\pi'}^{\star})/ B^k(\mathcal{A}_{\pi'}^{\star}) 
 \quad \text{and} \quad
 \overline{H}^k(\mathcal{A}_{\pi'}^{\star}) = Z^k(\mathcal{A}_{\pi'}^{\star})/ \overline{B}^k(\mathcal{A}_{\pi'}^{\star}) .
$$

 \medskip
 
\begin{proposition}\label{pp.duality1}
 The pairing given by (\ref{eq:pairing}) induces a duality between the reduced cohomology
of \(\mathcal{A}_{\pi'}^{\star}\) and that of  \({\Sigma}_{\pi}^{\star}(M)\):
$$
 \overline{H}^{k-1}(\mathcal{A}_{\pi'}^{\star}) \times \overline{H}^{n-k}({\Sigma}_{\pi}^{\star}(M))
 \to \r.
$$
\end{proposition}

\textbf{Proof.}
Observe first that
\begin{equation}\label{eq.orths1}
   B^{k-1}(\mathcal{A}_{\pi'}^{\star}) \subset Z^{k-1}(\mathcal{A}_{\pi'}^{\star})= (B^{n-k}(\mathcal{A}_{\pi'}^{\star})])^\perp \subset \mathcal{A}_{\pi'}^{k-1}.
\end{equation} 
Indeed, both inclusions are trivial and the above  equality follows from the fact that   $d_{\Sigma}$ and $d_{\mathcal{A}}$ are adjoint operators. 

Let us recall the (classic) argument. 
Fix  \(a \in \mathcal{A}_{\pi'}^{k-1}\), if  $a \in \ker[d_{\mathcal{A}} ^{k-1}]$, then $\langle a, d_{\Sigma}s\rangle = \pm \langle d_{\mathcal{A}} a, s\rangle = 0$ for any   \(s\in {\Sigma}_{\pi}^{n-k-1}(M)\), i.e. $a \in  (\im[d_{\Sigma}^{n-k-1}])^\perp$. Conversely, if  $a \in  (\im[d_{\Sigma}^{n-k-1}])^\perp$, then 
$ \langle d_{\mathcal{A}} a, s\rangle = \pm \langle a, d_{\Sigma}s\rangle = 0$ for any $s\in {\Sigma}_{\pi}^{n-k-1}(M)$
because  \({\Sigma}_{\pi}^{n-k-1}(M)\)   is dual to \(\mathcal{A}_{\pi'}^{k}\). This means that $d_{\mathcal{A}} a =0$.
 
 \medskip
 
Similarly, by Theorem \ref{th.propertyofA}  we have
\begin{equation}\label{eq.orths2}
  \im[d_{\Sigma}^{n-k-1}] \subset  \ker[d_{\Sigma}^{n-k}]= (\im[d_{\mathcal{A}} ^{k-2}])^\perp 
  \subset {\Sigma}_{\pi}^{n-k}(M).
\end{equation} 

The Proposition follows now from Lemma \ref{lem.dualquotient} for $X_0=\mathcal{A}_{\pi'}^{k-1}$, $X_1={\Sigma}_{\pi}^{n-k}(M)$ and equations (\ref{eq.orths1}) and (\ref{eq.orths2}) since by definition we
have
$$
\overline{H}^{k-1}(\mathcal{A}_{\pi'}^{\star}) = \ker[d_{\mathcal{A}} ^{k-1}]/\overline{\im[d_{\mathcal{A}} ^{k-2}]}
$$
and
$$
\overline{H}^{n-k}({\Sigma}_{\pi}^{\star}(M)) = \ker[d_{\Sigma}^{n-k}]/\overline{\im[d_{\Sigma}^{n-k-1}]}.
$$
\qed
 
 \bigskip

\begin{proposition}  \label{prop.isomAA}
The reduced and non reduced cohomology of the Banach complex 
\(\left(\mathcal{A}_{\pi'}^{\star} , d_{\mathcal{A}}\right)\) is isomorphic to the interior 
$L_{\pi'}$-cohomology of $M$ up to a shift:
\begin{equation}\label{eq.isomAA}
  H^k_{\pi';0}(M) = H^{k-1}(\mathcal{A}_{\pi'}^{\star}) 
  \quad and \quad
  \overline{H}^k_{\pi';0}(M) = \overline{H}^{k-1}(\mathcal{A}_{\pi'}^{\star}) . 
\end{equation}
These isomorphisms  are induced from the map   $j : Z^{k}_{\pi';0}(M) \to \mathcal{P}_{\pi'}^{k-1}$
defined by $j(\beta) = \binom{0}{\beta}$ .
\end{proposition}

Because we have a short exact sequence of complexes
$$
 0 \to \Omega_{\pi';0}^{\star}(M)  \to \mathcal{P}_{\pi'}^{\star}  \to  \mathcal{A}_{\pi'}^{\star} \to 0
$$
and $\mathcal{P}_{\pi'}^{\star}$ has trivial cohomology, the result follows from general principles (see e.g.  Theorem 1b in \cite{KS99}).
However we give below a more informative, explicit  proof. 

\medskip

\textbf{Proof.}
 It will be convenient to describe the cohomology $H^{k-1}(\mathcal{A}_{\pi'}^{\star})$ as a 
quotient of $\mathcal{P}_{\pi'}^{k-1}$. An element in $\mathcal{A}_{\pi'}^{k-1}$ is represented
by an element $\binom{\alpha}{\beta} \in \mathcal{P}_{\pi'}^{k-1}$ modulo ${\Sigma}^{k-1}_{\pi';0}$,
i.e. the equality 
$$
 \binom{\alpha}{\beta} = \binom{\alpha + \varphi}{\beta + d\varphi}
$$
holds in $\mathcal{A}_{\pi'}^{k-1}$ if and only if $\varphi \in {\Sigma}^{k-1}_{\pi';0}$.

\smallskip

Now $\binom{\alpha}{\beta}$ represents an element in $Z^{k-1}(\mathcal{A}_{\pi'}^{\star})$
if $d_{\mathcal{A}}\binom{\alpha}{\beta} = \binom{\beta}{0}  \in {\Sigma}^{k}_{\pi';0}$.
This means that  $\beta  \in \Omega^{k}_{\pi';0}(M)$ and $d\beta = 0$, i.e.
$\beta  \in Z^{k}_{\pi';0}(M)$, in other words we have
$$
 Z^{k-1}(\mathcal{A}_{\pi'}^{\star}) = \left\{\binom{\alpha}{\beta} \in \mathcal{P}_{\pi'}^{k-1} \tq
  \beta  \in Z^{k}_{\pi';0}(M) \right\}  \Big{/}  {\Sigma}^{k-1}_{\pi';0}.
$$
Likewise, $\binom{\alpha}{\beta}$ represents an element in $B^{k-1}(\mathcal{A}_{\pi'}^{\star})$
if there exists $\binom{\gamma}{\delta} \in \mathcal{P}_{\pi'}^{k-2}$ such that 
$\binom{\alpha}{\beta} = d_{\mathcal{A}} \binom{\gamma}{\delta}=
 \binom{\delta}{0}$ modulo ${\Sigma}^{k-1}_{\pi';0}$, and this means that 
 $\beta = d\varphi  \in B^{k}_{\pi';0}(M)$. We thus have
$$
 B^{k-1}(\mathcal{A}_{\pi'}^{\star}) = \left\{\binom{\alpha}{\beta} \in \mathcal{P}_{\pi'}^{k-1} \tq
  \beta  \in B^{k}_{\pi';0}(M) \right\}  \Big{/}  {\Sigma}^{k-1}_{\pi';0},
$$
and, by completion
$$
 \overline{B}^{k-1}(\mathcal{A}_{\pi'}^{\star}) = \left\{\binom{\alpha}{\beta} \in \mathcal{P}_{\pi'}^{k-1} \tq
  \beta  \in \overline{B}^{k}_{\pi';0}(M) \right\}  \Big{/}  {\Sigma}^{k-1}_{\pi';0}.
$$
The above equalities allows us to give the following explicit description of the cohomology
of $\mathcal{A}_{\pi'}^{\star}$:
$$
 H^{k-1}(\mathcal{A}_{\pi'}^{\star}) =  \left\{\binom{\alpha}{\beta} \in \mathcal{P}_{\pi'}^{k-1} \tq
  \beta  \in Z^{k}_{\pi';0}(M) \right\} 
  \Big{/}  \left\{\binom{\lambda}{\mu} \in \mathcal{P}_{\pi'}^{k-1} \tq
  \mu \in B^{k}_{\pi';0}(M) \right\}.
$$
But this quotient can also be described as
$$
 H^{k-1}(\mathcal{A}_{\pi'}^{\star}) =  \left\{\binom{0}{\beta} \in \mathcal{P}_{\pi'}^{k-1} \tq
  \beta  \in Z^{k}_{\pi';0}(M) \right\} 
  \Big{/}  \left\{\binom{0}{\mu} \in \mathcal{P}_{\pi'}^{k-1} \tq
  \mu \in B^{k}_{\pi';0}(M) \right\}.
$$
In short, we have established that the embedding $j : Z^{k}_{\pi';0}(M) \to \mathcal{P}_{\pi'}^{k-1}$
defined by $j(\beta) = \binom{0}{\beta}$ induces an algebraic isomorphism
$j : H^k_{\pi';0}(M) \cong  H^{k-1}(\mathcal{A}_{\pi'}^{\star})$.
We also have 
$$
 \overline{H}^{k-1}(\mathcal{A}_{\pi'}^{\star}) =  \left\{\binom{0}{\beta} \in \mathcal{P}_{\pi'}^{k-1} \tq
  \beta  \in Z^{k}_{\pi';0}(M) \right\} 
  \Big{/}  \left\{\binom{0}{\mu} \in \mathcal{P}_{\pi'}^{k-1} \tq
  \mu \in \overline{B}^{k}_{\pi';0}(M) \right\}.
$$
This quotient is equipped with its natural quotient norm and the homomorphism $j$ clearly induces
an isometric isomorphism 
$j : \overline{H}^k_{\pi';0}(M) \cong  \overline{H}^{k-1}(\mathcal{A}_{\pi'}^{\star})$.

\qed

\medskip

From   the propositions \ref{prop.isomAA} and  \ref{pp.duality1}, we now deduce 
the following duality result:

\begin{theorem} \label{thm:th.HPduality} 
Let $M$ be an arbitrary smooth $n$-dimensional oriented Riemannian manifold.
For any $1<\pi<\infty$,  the Banach spaces $\overline{H}_{\pi}^{k}(M)$ and $\overline{H}_{\pi';0}^{n-k}(M)$ 
are in duality for the pairing $\left\langle \beta,\omega \right\rangle = \int_M  \beta\wedge \omega.$
where $\beta\in Z_{\pi}^{k}(M)$ and $\omega \in Z_{\pi';0}^{n-k}(M)$. \\
\end{theorem}

\textbf{Proof.} By Propositions  \ref{pp.duality1}, we know that $\overline{H}^{k-1}(\mathcal{A}_{\pi'}^{\star})$ is isomorphic
to the dual of $\overline{H}^{n-k}({\Sigma}_{\pi}^{\star}(M))$ for the pairing  (\ref{eq:pairing}):
$$
\left\langle \binom{\alpha}{\beta} , \binom{\omega}{d\omega}  \right\rangle =
\int_M ((-1)^{k}\alpha\wedge d\omega + \beta\wedge \omega).
$$
By Propositions \ref{prop.isomAA}, we have an isomorphism $\overline{H}^k_{\pi';0}(M) \cong \overline{H}^{k-1}(\mathcal{A}_{\pi'}^{\star})$ induced by the map  $\beta \mapsto \binom{0}{\beta}$, and we trivially have an isomorphism
$\overline{H}_{\pi';0}^{n-k}(M) \cong \overline{H}^{n-k}({\Sigma}_{\pi}^{\star}(M))$ given by $\omega \mapsto \binom{\omega}{0}$.
It follows that  $\overline{H}_{\pi}^{k}(M)$ and $\overline{H}_{\pi';0}^{n-k}(M)$ 
are in duality for the pairing 
$$
\left\langle \binom{0}{\beta} , \binom{\omega}{0}  \right\rangle =
\int_M  \beta\wedge \omega.
$$
\qed

\medskip

Observe that Theorem \ref{th.duality1} is now also proved, since it suffices to apply the previous Theorem to any
sequence $\pi=\{ p_{0},p_{1},\cdots,p_{n}\}\subset (1,\infty)$ such that $q=p_{k-1}$ and $p = p_k$.

\section{Parabolicity}

In this section, we prove Theorem \ref{th.parcomp}. Recall the definition of an  \(s\)-\emph{parabolic} manifold.

\medskip

\begin{definition}  
The Riemannian manifold  \((M,g)\) is said to be  \(s\)-\emph{parabolic}, \(1\leq s \leq \infty\)  if there exists a sequence of smooth functions with compact support \(\{\eta_j\} \subset C^{\infty}_0(M)\) such that 
\begin{enumerate}[i.)]
  \item \(0 \leq \eta_j \leq 1\);
  \item \(\eta_j \to 1\) uniformly on every compact subset of \(M\);
  \item \(\ds \lim_{j\to \infty} \int_M |d\eta_j|^s = 0\).
\end{enumerate}
A non \(s\)-parabolic manifold is called \(s\)-\emph{hyperbolic}.
\end{definition} 

It is easy to check that $M$ is $\infty$-parabolic if and only
if it is complete. If the manifold is complete with finite volume,
then it is $s$-parabolic for any $s\in[0,\infty]$. There are numerous
characterizations of $s$-parabolicity as well as a number of geometric
interpretations. See \cite{Holop99,tr99,GoT99} and references therein for more this subject.

\bigskip

\begin{proposition}\label{prop.parab1}
If \(M\) is \(s\)-parabolic and  \(\frac{1}{s} =  \frac{1}{p} - \frac{1}{q}\), then $\mathcal{D}^{\ell}$ is dense in \(\Omega_{q,p}^{\ell}(M)\).
\end{proposition}

\textbf{Proof.}
Let $\alpha\in\Omega_{q,p}^{\ell}(M)$. We know by Lemma \ref{lem.density1}
that there exists a sequence $\alpha_{j}\in C^{\infty}(M,\Lambda^{\ell})$
converging to $\alpha$ in $\Omega_{q,p}^{\ell}(M)$. Since $M$ is
$s$-parabolic, there exists a sequence $\eta_{j}\in C_{0}^{\infty}(M)$
satisfying the conditions (i)--(iii) above. It is then clear that
$\beta_{j}=\eta_{j}\alpha_{j}\in \mathcal{D}^{\ell}$,
and we claim that $\beta_{k}$ converges to $\alpha$ in $\Omega_{q,p}^{\ell}(M)$.
Indeed, we have 
\begin{align*}
\Vert\alpha-\beta_{j}\Vert_{L^{q}} & \leq\Vert\alpha-\alpha_{j}\Vert_{L^{q}}+\Vert\alpha_{j}-\beta_{j}\Vert_{L^{q}}\\
 & =\Vert\alpha-\alpha_{j}\Vert_{L^{q}}+\Vert(1-\eta_{j})\alpha_{j}\Vert_{L^{q}}.
\end{align*}
Since $\Vert\alpha-\alpha_{j}\Vert_{L^{q}}\rightarrow0$ by hypothesis
and $\Vert(1-\eta_{j})\alpha_{j}\Vert_{L^{q}}\rightarrow0$ by the
Lebesgue Dominated convergence theorem, we have $\Vert\alpha-\beta_{j}\Vert_{L^{q}}\rightarrow0$.
We also have 
\begin{align*}
\Vert d\alpha-d\beta_{j}\Vert_{L^{p}} & \leq\Vert d\alpha-d\alpha_{j}\Vert_{L^{p}}+\Vert d\alpha_{j}-d\beta_{j}\Vert_{L^{p}}\\
 & =\Vert d\alpha-d\alpha_{j}\Vert_{L^{p}}+\Vert d\left((1-\eta_{j})\alpha_{j}\right)\Vert_{L^{p}}\\
 & =\Vert d\alpha-d\alpha_{j}\Vert_{L^{p}}+\Vert(1-\eta_{j})d\alpha_{j}\Vert_{L^{p}}+\Vert d(1-\eta_{j})\wedge\alpha_{j}\Vert_{L^{p}}\\
 & =\Vert d\alpha-d\alpha_{j}\Vert_{L^{p}}+\Vert(1-\eta_{j})d\alpha_{j}\Vert_{L^{p}}+\Vert d(1-\eta_{j})\Vert_{L^{s}}\Vert\alpha_{j}\Vert_{L^{q}}.
\end{align*}
Which converges to $0$. We have shown that for any $\alpha\in\Omega_{q,p}^{\ell}(M)$,
there exists a sequence $\beta_{j}\in\mathcal{D}^{\ell}$
such that 
\[
\lim_{j\rightarrow0}\Vert\alpha-\beta_{j}\Vert_{q,p}=\Vert\alpha-\beta_{j}\Vert_{L^{q}}+\Vert d\alpha-d\beta_{j}\Vert_{L^{p}}=0.
\]
 This means that $\mathcal{D}^{\ell}\subset \Omega_{q,p}^{\ell}(M)$ is dense. 
 
\qed

\bigskip

\textbf{Proof of  Theorem \ref{th.parcomp}.}
The previous Proposition implies that if  \(M\) is \(s\)-parabolic for \(\frac{1}{s} =  \frac{1}{p} - \frac{1}{q}\), 
then the set of smooth exact $k$-forms with compact support is dense in $B^k_{q,p}(M)$. In particular, $\overline{B}^k_{q,p}(M)$ is the closure of $d\mathcal{D}^{k-1}(M)$ in $L^p(M,\Lambda^k)$:
$$
 \overline{B}^k_{q,p}(M) = \overline{d\mathcal{D}^{k-1}(M)\ }^{L^p}.
$$
Therefore
$$
  \overline{H}^k_{q,p}(M) = Z^k_{p}(M)/\overline{B}^k_{q,p}(M) =
  Z^k_{p}(M)/\overline{d\mathcal{D}^{k-1}(M)\ } = \overline{H}^{k}_{\comp,p}(M).
$$
The  Proposition (\ref{prop.parab1}) also implies that any closed  $k$-form in  
$Z^k_{q}(M)$ can be approximated in the $\Omega_{q,p}^{k}(M)$-topology
by a sequence of smooth  $k$-forms with compact support. 
Hence 
$$
 Z^k_{q}(M) = \ker(d)\cap \overline{\mathcal{D}^{k}(M)\ }^{\Omega_{q,p}^{k}} =  Z^k_{q,p;0}(M),
$$
and we have
$$
  \overline{H}^k_{q,p;0}(M) = Z^k_{q,p;0}(M)/\overline{B}^k_{q;0}(M) =
  Z^k_{q}(M)/\overline{d\mathcal{D}^{k-1}(M)\ } = \overline{H}^{k}_{\comp,q}(M).
$$
The Theorem is proved.

\qed
 
\section{The conformal cohomology}

\textbf{Definition} The \emph{conformal de Rham complex} of the $n$-dimensional 
 Riemannian manifold $(M,g)$ is defined as
$$
  \Omega_{\conf}^{k}(M) = \{ \omega \in L^{n/k}(M,\Lambda^k) \tq d\omega \in L^{n/(k+1)}(M,\Lambda^{k+1}) \}.
$$
It is thus simply the complex  $\Omega_{\pi}^{k}(M)$ associated to the sequence $\pi$ defined by $p_k = n/k$.
It is proved in \cite{GT2009a} that $\Omega_{\conf}^{k}(M)$ is a differential graded Banach algebra and is invariant 
under quasi-conformal maps. Moreover, a homeomorphism $f : (M,g) \to (N,h)$  between two Riemannian manifolds.  is  a quasiconformal map if and only if the pull-back of differential forms defines an isomorphism of Banach 
differential algebras
 $
  f^*  : \Omega_{\conf}^{\bullet}(N) {\longrightarrow} \Omega_{\conf}^{\bullet}(M).
 $
 
 \medskip
 
The conformal sequence $\pi = (p_k)$ defined by $p_k  = n/k$ is its own dual, because 
$\frac{1}{p_k}+\frac{1}{p_{n-k}} = \frac{k}{n}+\frac{n-k}{k} = 1$. 
Theorem (\ref{th.duality1}) can thus be rested as follow in the particular case of the conformal cohomology:
\begin{theorem} 
Let $M$ be an  oriented  Riemannian manifold of dimension $n$ and  $2 \leq k \leq (n-2)$. 
Then  $\overline{H}^{k}_{\conf}(M)$ is isomorphic to the dual of  
$\,Ê\overline{H}^{n-k}_{\conf;0}(M)$
The duality is induced by the integration pairing (\ref{eq.chmpairing}).
\end{theorem}

\section{Applications}
In this section, we give some consequences of the previously established results.

\begin{proposition}
 If  the  oriented Riemannian $(M,g)$  admits an isometry $f : M \to M$ without any bounded orbit,
 and $1<q,p<\infty$,  then  $\overline{H}^k_{q,p}(M)$  is finite dimensional if and only if 
 $\overline{H}^{k}_{q,p}(M) = 0$. The same holds for $\overline{H}^{k}_{\comp,p}(M)$.
\end{proposition}  

Gromov made this observation for $L^p$-cohomology assuming $M$ to be complete and of 
bounded geometry, see \cite[p. 220]{Gromov93}.  

\medskip

\textbf{Proof.} 
We essentially follow Gromov's argument.  Assume 
 that $\overline{H}^k_{q,p}(M))$ has finite dimension $m$, then we also have $m = \dim (\overline{H}^{n-k}_{p',q';0}(M)) < \infty$ by Theorem \ref{th.duality1}. Let $\varphi_1, \cdots , \varphi_m \in Z^{n-k}_{p',q';0}(M)$ be a set of closed forms representing a basis of 
$\overline{H}^k_{p',q';0}(M)$. By the Lebesgue dominated convergence theorem, we have for any $\omega \in Z^k_{p}(M)$
$$
 \lim_{j \to\infty } \int_M  (f^j)^*\omega \wedge \varphi_i = 0
$$
for  $i=1, \cdots, m$. This implies that $\lim_{j \to\infty } [(f^j)^*\omega] = 0$ in $\overline{H}^k_{q,p}(M)$
but  $f$ is an isometry, therefore
$$
 \|[\omega]\| = \lim_{j \to\infty }\| [(f^j)^*\omega]\| = 0
$$
for any element $[\omega] \in \overline{H}^k_{q,p}(M)$.
The proof for $\overline{H}^{k}_{\comp,p}(M)$ is the same.

\qed

\medskip

Our next result gives us an explicit criterion implying that a given cohomology class is non zero.

\begin{proposition}
Let \((M,g)\) be an oriented Riemannian manifold of dimension \(n\) and let  \ \(\alpha \in Z_{p}^{k}(M)\).  
Then the following conditions are equivalent.
\begin{enumerate}[(a)]
  \item  \([\alpha ]\neq 0\) in \(\overline{H}_{q,p}^{k}(M)\).
  \item There exists  $\omega \in Z^{n-k}_{p',q';0}(M)$ such that 
\begin{equation} \label{nnpair}
 \int_{M} \alpha \wedge\omega \neq 0.
\end{equation}
  \item There exists a sequence \(\ \{\gamma _{i}\}\subset
\mathcal{D}^{n-k}\) such that \
\begin{enumerate}[i)]
 \item   \(\ds \limsup_{i\rightarrow \infty }  \int_{M}\alpha
\wedge \gamma_{i}>0\).
\item  \(\ds \lim_{i\rightarrow \infty }\left\| d\gamma _{i}\right\|_{q'}=0\) \
where \(q'=\frac{q}{q-1}\).
\item \(\left\| \gamma _{i}\right\|_{p'}\) is a
bounded sequence for  \(p'=\frac{p}{p-1}\).
\end{enumerate}
\end{enumerate}
\end{proposition}

\textbf{Proof.}
(a) $\Rightarrow$ (b).
If \([\alpha ]\neq 0\) in \(\overline{H}_{q,p}^{k}(M)\), then by  
the duality Theorem \ref{th.duality1}, there exists $[\omega]
\in  \overline{H}^{n-k}_{p',q';0}(M)$ such that (\ref{nnpair}) holds. The
cohomology class $[\omega]$ is represented by an element $\omega \in Z^{n-k}_{p',q';0}(M)$,
which is the desired differential form.

\medskip

(b) $\Rightarrow$ (c)
Let  $\omega \in Z^{n-k}_{p',q';0}(M)$ be an element such that  (\ref{nnpair}) holds. By 
 definition of $Z^k_{p',q';0}(M)$, there exists a sequence $\gamma_j \in \mathcal{D}^{k}(M)$ such that 
$\| \gamma_j - \omega\|_{p'} \to 0$ and $\| d\gamma_j \|_{q'} \to 0$. This sequence clearly satisfies the conditions
(i), (ii) and (iii).

\medskip

The implication (c) $\Rightarrow$ (a) has been proved in  \cite[sec. 8]{GT2006}. We repeat the proof for convenience.  
Suppose that $\alpha\in\overline{B}_{q,p}^{k}(M)$.
Then ${\displaystyle \alpha=\lim_{j\rightarrow\infty}d\beta_{j}}$
for $\beta_{j}\in L^{q}(M,\Lambda^{k-1})$ with $d\beta_{j}\in L^{p}(M,\Lambda^{k})$.
We have for any $i,j$ 
\[
\int_{M}\gamma_{i}\wedge\alpha=\int_{M}\gamma_{i}\wedge d\beta_{j}+\int_{M}\gamma_{i}\wedge(\alpha-d\beta_{j})\,\,.
\]
 For each $j\in\Bbb{N}$, we can find $i=i(j)$ large enough so that
\ $\Vert d\gamma_{i(j)}\Vert_{{q}'}\,\Vert\beta_{j}\Vert_{q}\leq1/j$,
we thus have 
\[
\quad\left|\int_{M}\gamma_{i(j)}\wedge d\beta_{j}\right|\quad\leq\quad\left|\int_{M}d\gamma_{i(j)}\wedge\beta_{j}\right|\quad\leq\Vert d\gamma_{i(j)}\Vert_{{q}'}\,\Vert\beta_{j}\Vert_{q}\leq\frac{1}{j}\,.
\]
 On the other hand 
 \[
\quad\lim_{j\rightarrow\infty}\left|\int_{M}\gamma_{i(j)}\wedge(\alpha-d\beta_{j})\right|\quad
\leq\lim_{j\rightarrow\infty}\Vert\gamma_{i(j)}\Vert_{{p}^{\prime}}\,\Vert(\alpha-d\beta_{j})\Vert_{p}=0\,
\]
since \ $\Vert\gamma_{i(j)}\Vert_{{p}'}\,$\ is a bounded sequence
and $\,\,\Vert(\alpha-d\beta_{j})\Vert_{p}\rightarrow0$. It follows
that $\int_{M}\gamma_{i(j)}\wedge\alpha\rightarrow0$ in contradiction
to the hypothesis.

\qed
 
\bigskip 

The previous criterion is not always very practical because we cannot always produce 
useful closed forms with compact support. The next result is thus convenient in the case
of complete manifolds.

\medskip
 
\begin{corollary}
Assume that \(M\) is complete. Let \(\alpha\in
Z_{p}^{k}(M)\), and assume that there exists a  closed
\((n-k)\)-form \(\gamma\in Z_{p'}^{n-k}(M)\cap Z_{q'}^{n-k}(M)\),
where \(p'=\frac{p}{p-1}\) and \(q'=\frac{q}{q-1}\), such that
\[
\int_{M}\gamma\wedge\alpha > 0,
\]
then \(\alpha\,\notin\,\overline{B}_{q,p}^{k}(M)\) where
\(q'=\frac{q}{q-1}\). In particular,
\(\overline{H}_{q,p}^{k}(M)\neq0\).
\end{corollary}

This result is  Proposition 8.4 from \cite{GT2006}, an example of its 
usefulness can be found in that paper where it is used to prove the non vanishing of the 
$L_{q,p}$-cohomology of the hyperbolic plane. We give a short proof.

\medskip

\textbf{Proof.}   Let $\gamma$ be a differential form
as in the statment, by density we may assume $\gamma$ to be smooth.
Let $\{\eta_j\} \subset C^{\infty}_0(M)$ be a sequence of smooth function with compact support uniformly converging to $1$
and such that $\| \eta_j\|_{L^{\infty}} \to 0$. Then the sequence  $\gamma_j = \eta_j\ \gamma$ satisfies the condition (c) of the Proposition.

\qed

\medskip

If $p=q$, then the previous condition is not only sufficient, but also necessary.

\medskip 
 
\begin{corollary}
 Let $M$ be a complete oriented Riemannian manifold. Then 
an element  \(\alpha \in Z_{p}^{k}(M)\) is non zero in reduced $L_p$-cohomology 
if and only if there exists  $\omega \in Z^{n-k}_{p'}(M)$ such that 
$$
 \int_{M} \alpha \wedge\omega \neq 0.
$$
\end{corollary}

\textbf{Proof.} This is a special case of the previous Proposition, since
for complete manifolds, we have $Z^{n-k}_{p'}(M) = Z^{n-k}_{p',p';0}(M)$.
\qed
 
 \medskip
 
Based on the latter Corollary, we can prove the following vanishing result which has been 
first observed by Gromov, see also \cite[Proposition 15]{Pansu2008}.

\medskip

\begin{proposition}
 Assume that the Riemannian manifold is complete and admits a complete and proper Killing vector field $\xi$, then $\overline{H}^{k}_{\comp,p}(M)=0$  and (therefore) $\overline{H}^{k}_{q,p}(M)=0$ for any $1<p<\infty$ and $1 \leq q < \infty$.
\end{proposition}

\medskip

For instance if $N$ is a complete Riemannian manifold, then  $\overline{H}^{k}_{q,p}(N\times \r)=0$ for any $1<p<\infty$ and $1 \leq q <  \infty$.

\medskip
 
\textbf{Remark.} The hypothesis $p \neq 1$ is a necessary condition. For instance 
$\overline{H}^{k}_{q,1}(\r) \neq 0$ for any $1 \leq q < \infty$ (see \cite[prop. 9.3]{GT2006}).

\medskip

\textbf{Proof.} 
Because $M$ is complete, we have  $\overline{H}^{k}_{\comp,p}(M) =  \overline{H}^{k}_{p}(M)$. We thus only
need to prove the vanishing of the $L_p$-cohomology of $M$. The proof goes as follow: Because $\xi$ is a Killing vector field, we have $\|i_{\xi} \alpha \|_{L^p} \leq \|\alpha \|_{L^p}$ for any differential
form $\alpha$. Let   $f_t$ be the flow of  $\xi$. Using Cartan's formula $\mathcal{L}_{\xi}  \alpha  = di_{\xi} \alpha + i_{\xi} d\alpha$, we see that for any 
$\alpha \in Z^{k}_{p}(M)$ and  $\varphi \in Z^{n-k}_{p'}(M)$ we have
$$
 \frac{d}{dt} \int_M f_t^* \alpha \wedge \varphi = \int_M \mathcal{L}_{\xi}  \alpha  \wedge \varphi 
 = \int_M  d i_{\xi}  \alpha  \wedge \varphi  = 0
$$
because $[d i_{\xi}  \alpha ] = 0$ in $\overline{H}^{k-1}_{p}(M)$. On the other hand, because the flow $f_t$ is proper, we have
$$
 \lim_{t\to 0} \int_M f_t^* \alpha \wedge \varphi   = 0.
$$
It follows that $ \int_M f_t^* \alpha \wedge \varphi   = 0$ for any $\alpha \in Z^{k}_{p}(M)$ and  $\varphi \in Z^{n-k}_{p'}(M)$
and therefore $\overline{H}^{k}_{p}(M)= 0$.
\qed


\end{document}